\documentclass[preprint,11pt]{elsarticle}
\usepackage[table]{xcolor}
\usepackage{comment}
\usepackage{tikz}
\usetikzlibrary{shapes,arrows,positioning}
\usetikzlibrary{arrows.meta}
\usepackage[margin=1in]{geometry}
\usepackage{amssymb,amsmath}
\usepackage{amsthm}
\usepackage{bm}
\usepackage{setspace}
\usepackage{enumitem}
\usepackage{graphicx}
\usepackage{url}
\usepackage{booktabs}
\usepackage{csvsimple}
\usepackage{multicol,multirow}
\usepackage{subcaption}
\usepackage[export]{adjustbox}
\usepackage[scaled=0.9]{newpxtext}
\usepackage[scaled=0.9]{newpxmath}

\usepackage{hyperref}
\hypersetup{
  colorlinks =true,
  linkcolor = blue,
  filecolor = blue,
  citecolor = magenta,  
  urlcolor = magenta,
  }

\theoremstyle{remark}

\usepackage{subcaption}
\usepackage{float}
\floatstyle{ruled}
\newfloat{model}{thp}{lop}
\floatname{model}{Model}

\usepackage{multirow,mdframed}
\mdfsetup{skipabove=0pt,skipbelow=0pt}
\mdfdefinestyle{model}{%
  innertopmargin=0pt,
  innerbottommargin=0pt,
  innerleftmargin=0pt,
  innerrightmargin=0pt,
  skipbelow=0pt,%
  skipabove=0pt,
  splitbottomskip=0pt,
  splittopskip=0pt,
  leftmargin =0pt,%
    rightmargin=0pt,
  splittopskip=0pt,
    usetwoside=false,
}
\usepackage{setspace}
\allowdisplaybreaks 
\journal{Elsevier}
\begin{document}

\tikzstyle{decision} = [diamond, draw, fill=blue!20, 
text width=7em, text badly centered, inner sep=0pt]
\tikzstyle{block} = [rectangle, draw, fill=blue!20, 
text width=12em, text centered, rounded corners, minimum height=4em]
\tikzstyle{block2} = [rectangle, draw, fill=yellow!20, 
text width=12em, text centered, rounded corners, minimum height=4em]
\tikzstyle{line} = [draw, -latex']

\begin{frontmatter}
\title{Stability-Constrained AC Optimal Power Flow--A Gaussian Process-Based Approach}

\author[inst1]{Vincenzo~Di~Vito}
\affiliation[inst1]{organization={University of Virginia},
            addressline={Computer Science Department}, 
            city={Charlottesville},
            postcode={22903}, 
            state={VA},
            country={USA}}

\author[inst2]{Kaarthik~Sundar\corref{cor}}
\affiliation[inst2]{organization={Los Alamos National Laboratory},
            city={Los Alamos},
            postcode={87545}, 
            state={NM},
            country={USA}}
\author[inst1]{Ferdinando~Fioretto}
\author[inst3]{Deepjyoti~Deka}
\affiliation[inst3]{organization={Massachusetts Institute of Technology},
            addressline={MIT Energy Initiative}, 
            city={Cambridge},
            postcode={02139}, 
            state={MA},
            country={USA}}


\begin{abstract}
The Alternating Current Optimal Power Flow (ACOPF) problem is a core task in power system operations, aimed at determining cost-effective generation dispatch while satisfying physical and operational constraints. However, conventional ACOPF formulations rely on steady-state models and neglect generator dynamics, which can result in operating points that are economically optimal but dynamically unstable. This paper proposes a novel, data-driven approach to incorporate generator dynamics into the ACOPF using Gaussian Process (GP) models. Specifically, it introduces an exponential surrogate function to characterize the stability of solutions to the differential equations governing synchronous generator dynamics. The exponent, which indicates whether system trajectories decay (stable) or grow (unstable), is learned as a function of the bus voltage using GP regression. Crucially, the framework enables probabilistic stability assessment to be integrated directly into the optimization process. The resulting dynamics-aware ACOPF formulation identifies operating points that satisfy both operational safety and dynamic stability criteria. Numerical experiments on the IEEE 39-bus, 57-bus, and 118-bus systems demonstrate that, compared with existing data-driven approaches, the proposed method efficiently captures generator dynamics with limited training data, yielding more reliable and robust decisions across a wide range of operating conditions.
\end{abstract}




\end{frontmatter}

\section{Introduction} \label{sec:intro}
The AC Optimal Power Flow (ACOPF) problem in power system operations is formulated to determine economically optimal generation setpoints that satisfy system demand while adhering to steady-state physical laws, such as Kirchhoff’s current and voltage laws, and operational constraints, such as thermal and generator limits \cite{Vittalbook}. Despite its widespread use, the ACOPF formulation inherently assumes a static system and neglects the dynamic behavior of generators and system-wide frequency response. This abstraction, while computationally tractable and useful for routine operational planning, fails to account for transient and dynamic performance, which are critical for maintaining system stability. Recent work~\cite{divito2024learningoptimizemeetsneuralode} has demonstrated that learning-based approaches relying solely on steady-state representations, without accounting for system dynamics, can result in operating points that violate stability requirements. In practice, when the dispatch solutions from ACOPF lead to dynamically unstable operating points, system operators resort to heuristic, experience-based adjustments to stabilize the grid. These ad hoc interventions introduce uncertainty, elevate the risk of load curtailment, and often result in higher operational costs. Such issues are further amplified by the increasing penetration of renewable energy sources, which introduce significant variability and uncertainty~\cite{9286772}. 

The gap between steady-state optimization and dynamic system requirements highlights the need for integrated, stability-aware optimization frameworks in both operational practice and regulatory policy. In this context, it is important to distinguish between two key notions of stability relevant to system operators \cite{kundurbook}: (a) generator rotor angle stability, which is local to individual generators, and (b) small-signal stability, which concerns the overall system's dynamic response to small perturbations. While incorporating both is essential for comprehensive dynamic security, this work focuses on the former, with small-signal stability left for future investigation.

\textcolor{black}{However, incorporating model dynamics into the ACOPF formulation poses significant computational and modeling challenges. In particular, (1) the increased problem dimensionality and complexity of introducing differential equations within an already complex nonlinear, nonconvex optimization framework, and (2) the need for surrogate models that can accurately capture generator dynamics across varying operating conditions, while remaining computationally efficient. }
To address these challenges, this paper proposes to use Gaussian Process (GP) models as surrogates for generator dynamics. GPs offer two key advantages: (i) they provide a non-parametric continuous differentiable framework capable of modeling nonlinear behaviors and (ii) they provide principled uncertainty quantification of predicted variables. These features enable a natural integration of GP into ACOPF via chance constraints, thereby enabling probabilistic enforcement of each generator's rotor-angle stability criteria. 
In the context of modeling system dynamics, GPs have been employed in the literature to learn solutions of specific instances of ordinary differential equations (ODEs), typically under fixed system parameters \cite{heinonen2018learningunknownodemodels}. Learning system dynamics that remain valid across a range of parameters using GPs is challenging, as the surrogate models must generalize over diverse system states, including variations in load levels and network topology, without sacrificing accuracy or computational tractability.

To enable such generalizability and computational efficiency, we shift our focus from modeling full dynamic trajectories directly to characterizing stability using an exponential surrogate function derived from the underlying differential equations governing generator dynamics. The exponent in this surrogate function reflects the growth or decay rate of the generator's trajectory and hence indicates whether the generator's operating point is dynamically stable or unstable. This exponent, modeled as a GP, is learned as a function of the local bus voltage and phase angles and succinctly captures the nonlinear relationship between system states and generator stability.   
Our results in Section \ref{subsec:predictive-performance} show that the aforementioned approach predicts the stability of trajectories with an accuracy of around $85$-$90$\% across a wide range of generator parameters, corroborating its effectiveness in capturing the dominant behavior of generator dynamics.
Furthermore, it enables efficient integration of dynamic stability constraints into the ACOPF, bridging the gap between steady-state optimization and dynamic system performance. Specifically for the $118$-IEEE bus system, the solutions of the standard ACOPF without accounting for system dynamics, yield unstable setpoints in $41\%$ of the loading scenarios, whereas our GP-based stability integration in ACOPF produces generator setpoints that satisfy the stability requirements in $85\%$ of the cases.

In summary, the main contributions of this article are as follows:
\begin{enumerate}
    \item GP surrogate for generator rotor angle stability -- We introduce an exponential stability indicator that summarizes the behavior of synchronous generator dynamics, learned efficiently via GP regression on local bus voltage data. This model captures nonlinear dependencies and enables the probabilistic enforcement of stability constraints into a steady-state ACOPF problem formulation. 
    \item GP-based stability-constrained ACOPF -- We build on the framework of \cite{divito2024learningoptimizemeetsneuralode}, which introduces the stability-constrained ACOPF model and proposes a purely learning-based approach to capture generator dynamics. In contrast, we present an approach that directly incorporates the stability constraints learned by a GP into the steady-state ACOPF and solves the resulting problem using off-the-shelf primal-dual interior-point solvers. 
    \item Empirical validation -- We validate our approach on the IEEE 39-, 57-, and 118-bus systems. Empirical results show that the classical ACOPF often yields unstable setpoints, whereas our method consistently produces dynamically stable dispatch solutions.
\end{enumerate}

The remainder of the paper is organized as follows. In Section~\ref{sec:lit-review}, we review prior work on GP-based surrogates for learning solutions of differential equations. We also provide an overview of extensions to the traditional ACOPF formulation that incorporate dynamic behavior into the optimization framework.
In Section \ref{sec:prob_formulation}, we present the mathematical formulation of the ACOPF problem, describe the generator dynamics, and introduce the stability-constrained ACOPF. In Section \ref{sec:gp-intro}, we offer a concise overview of GPs for supervised regression. This overview is not intended as a comprehensive primer on GPs, but rather as a brief introduction. In Section \ref{sec:method}, we first show that generator rotor-angle trajectories either damp out or blow up, a pattern well captured by a compact exponential envelope whose growth-rate parameter cleanly flags stability. We then train a GP to forecast the growth rate from operating conditions and incorporate this forecast into the ACOPF formulation as a chance constraint to enforce rotor angle stability. Finally, in Section \ref{sec:exp}, we present results of extensive computational experiments that corroborate the effectiveness of the proposed approach in computing operating points for the transmission system that satisfy rotor angle stability constraints of the generators in a probabilistic sense, followed by a conclusion and avenues for future work in Section \ref{sec:conclusion}. 

\section{Related work} \label{sec:lit-review}
We categorize the related literature into two primary areas: (1) the use of Gaussian Processes (GPs) for learning the solution or behavior of the solutions to differential equation systems, and (2) approaches for integrating stability constraints into AC Optimal Power Flow (ACOPF) formulations.

\subsection{GPs for Learning Solutions of Differential Equations}
GPs have become a powerful tool in system modeling due to their nonparametric nature, their ability to capture nonlinear dependencies, and their inherent uncertainty quantification~\cite{rasmussen2006gaussian}. In the context of dynamical systems, GPs have been employed to learn both explicit solution trajectories and latent system dynamics from data. Authors in \cite{raissi2017physics} pioneered the use of GPs for learning the solutions of differential equations directly from noisy observations. Similarly, the authors of \cite{heinonen2018learningunknownodemodels} introduced methods for inferring latent ordinary differential equation (ODE) models from data using GP priors, thereby enabling data-driven reconstruction of underlying physical processes.
Other works have enhanced GP-based learning with physical structure and constraints. Latent force models~\cite{8485787} combine mechanistic models with GP-based learning to capture hybrid physical and data-driven dynamics, while physics-informed GPs~\cite{hamelijnck2021spatio} integrate prior knowledge of system dynamics directly into the GP formulation to improve interpretability and extrapolation capability. For a comprehensive treatment of GP modeling and its applications to the control of dynamic systems, we refer to \cite{kocijan2016modelling}, which provides both theoretical insights and practical applications.
A common theme across these approaches is the use of GPs to avoid explicit numerical integration, which is often computationally intensive. In our work, we build on this literature by using GPs to learn a compact, stability-oriented surrogate---specifically, an exponential function whose rate parameter  captures the system's stability behavior. Unlike prior work that focuses on full trajectory reconstruction, we focus on extracting a stability metric sufficient for integration into a constrained optimization framework. 

We also remark that GPs have been applied to a wide range of tasks in power systems, including load forecasting~\cite{9435025}, learning steady-state power flow solutions~\cite{9552521} and optimization \cite{mitrovic2023data}, probabilistic estimation of power system dynamics~\cite{9813554, 10266787}, risk assessment~\cite{10614750,pareek2023fast}, and stability-based control~\cite{11003897}. We refer the reader to~\cite{tan2025gaussianprocessespowersystems} for a comprehensive survey on this topic. Next, we discuss how to review prior work on the integration of stability criteria within optimization frameworks in power grids.

\subsection{Integrating Stability Constraints into ACOPF}
The ACOPF problem is traditionally formulated as a nonlinear programming problem. A variety of solution techniques have been developed, ranging from interior-point methods~\cite{zimmerman2010matpower} and sequential quadratic programming, to more recent convex relaxation approaches~\cite{molzahn2019survey}. Over time, several steady-state extensions to the ACOPF formulation have been proposed to account for additional operational considerations. 
These include security-constrained OPF~\cite{capitanescu2011state}, contingency analysis~\cite{monticelli2007security}, and multi-period OPF~\cite{jabr2014robust}.  However, these extensions typically treat the power system as static, lacking integration of dynamic behavior. Prior research on incorporating stability constraints derived from dynamic behavior into ACOPF has predominantly relied on data-driven approaches. This reliance arises because formulating ACOPF with stability constraints expressed as Lyapunov conditions yields a nonlinear semi-definite optimization problem, for which, to the best of our knowledge, no open-source or commercial solvers are currently available. On the other hand, while methods for simulating dynamic trajectories for power grids and components exist, explicitly inserting the discretization techniques (e.g., collocation or finite‐difference schemes) used in such methods directly into the already complex ACOPF would result in prohibitively large and computationally intractable formulations. As such, surrogate models have been used to represent trajectories or stability criteria within the OPF formulation, and recent studies have proposed integrating data-driven system dynamics and stability constraints into OPF frameworks. Authors in \cite{liu2021explicit} introduced a small-signal stability-constrained OPF formulation that uses a support vector machine-based surrogate model to represent the stability boundary. Authors in \cite{su2023icnn} proposed an analytic optimization framework for transient stability-constrained OPF, where an Input Convex Neural Network is used as a surrogate model for transient stability assessment.  Most relevant to this work, the authors of \cite{divito2024learningoptimizemeetsneuralode} proposed a learning-based framework to solve a transient stability-constrained ACOPF problem that integrates synchronous generator and related stability constraints within a data-driven ACOPF model. 

Building on these advances, our work proposes a novel formulation that leverages efficient Gaussian Process surrogates to capture synchronous generator dynamics, for direct integration within the ACOPF model. These surrogates enable accurate enforcement of dynamic stability constraints within the ACOPF framework without requiring explicit trajectory simulation during optimization.

\section{Problem Formulation} \label{sec:prob_formulation}
This section introduces the stability-constrained ACOPF problem by first presenting two key components: the steady-state ACOPF formulation and the dynamics of synchronous generators. The paper adopts the following notation: lowercase letters denote scalars, whereas boldface letters denote vectors. Uppercase symbols represent complex variables, which may be expressed in either rectangular or polar form. Sets are denoted using standard calligraphic symbols.

\subsection{Steady-State ACOPF}
The steady-state ACOPF problem determines the most cost-effective generator dispatch that satisfies demand within a power network, subject to physical laws and operational limits. The power network is modeled as a graph $({\mathcal{N}, \mathcal{L}})$, where the node set $\mathcal{N}$ consists of $n$ buses and the edge set $\mathcal{L}$ represents $\ell$ transmission lines. The set of synchronous generators is denoted by $\mathcal{G} \subseteq \mathcal{N}$.

The AC power flow model employs complex-number representations of current $I$, voltage $V$, admittance $Y$, and power $S$ to capture both the physical behavior of power flow and system-level constraints. At each bus $i \in \mathcal{N}$, the generated and demanded power are denoted by $S^r_i = p_i^r + j q_i^r$ and $S^d_i = p_i^d + j q_i^d$, respectively. The power flow across line $(i, j) \in \mathcal{L}$ is represented as $S_{ij}$, and $\theta_i$ denotes the voltage phase angle at bus $i$.

The AC power flow physics is governed by Kirchhoff's and Ohm's laws: Kirchhoff's Current Law (KCL) is given by $I^r_i - I^d_i = \sum_{(i,j)\in \mathcal{L}} I_{ij}$, Ohm's Law is expressed as $I_{ij} = Y_{ij}(V_i - V_j)$, and the complex power flow is defined as $S_{ij} = V_i I_{ij}^*$. These relationships form the foundation of the AC power flow constraints, captured in equations~\eqref{eq:ac_3} and~\eqref{eq:ac_4} in Model~\ref{model:ac_opf}. The objective function~\eqref{ac_obj} minimizes the total generation cost, where each generator $i \in \mathcal{G}$ has a quadratic cost function defined by coefficients $c_{2i}, c_{1i}, c_{0i}$. Specifically, the cost is modeled as $c_{2i} (\Re(S^r_i))^2 + c_{1i} \Re(S^r_i) + c_{0i}$.

The model includes a comprehensive set of operational constraints. Voltage magnitude bounds are enforced by~\eqref{eq:ac_1}, and phase angle difference limits by~\eqref{eq:ac_6}. Generator output limits are defined in~\eqref{eq:ac_2}, and line flow limits are enforced via~\eqref{eq:ac_5}. Equation~\eqref{eq:ac_7} sets a reference for phase angle by fixing $\theta_{\text{ref}} = 0$. 
It is important to note that the formulation in Model~\ref{model:ac_opf} does \emph{not} include synchronous generator dynamics. While linearized, stability-aware OPF formulations have been proposed (e.g.,~\cite{Hafez16}), such dynamic considerations are absent in the standard steady-state ACOPF model.

\begin{model}[tbp]
    \caption{The Steady-State ACOPF Problem}
    \label{model:ac_opf}
    {\small
    \begin{subequations}
    \begin{align}
        \mbox{\bf variables:} \;\;
        & S^r_i, V_i \;\; \forall i\in \mathcal{N}, \;\;
          S_{ij}     \;\; \forall(i,j)\in \mathcal{L} \nonumber \\
        \mbox{\bf minimize:} \;\;
        & \sum_{i \in \mathcal{G}}  c_{2i} (\Re(S^r_i))^2 + c_{1i}\Re(S^r_i) +  c_{0i} \label{ac_obj} \\
        \mbox{\bf subject to:} \;\; 
        & v^l_i \leqslant |V_i| \leqslant v^u_i       \;\; \forall i \in \mathcal N \label{eq:ac_1} \\
        & -\theta^\Delta_{ij} \leqslant \angle (V_i V^*_j) \leqslant \theta^\Delta_{ij} \;\; \forall (i,j) \in \mathcal{L}  \label{eq:ac_6}  \\
        & S^{rl}_i \leqslant S^r_i \leqslant S^{ru}_i \;\; \forall i \in \mathcal{N} \label{eq:ac_2}  \\
        & |S_{ij}| \leqslant s^u_{ij}                  \;\; \forall (i,j) \in \mathcal{L} \label{eq:ac_5}  \\
        & S^r_i - S^d_i = \textstyle\sum_{(i,j)\in L} S_{ij} \;\; \forall i\in \mathcal{N} \label{eq:ac_3}  \\ 
        & S_{ij} =  Y^*_{ij} |V_i|^2 - Y^*_{ij} V_i V^*_j             \;\; \forall (i,j)\in \mathcal{L}
        \label{eq:ac_4} \\
        &\theta_{\text{ref}} = 0 \label{eq:ac_7}
    \end{align}
    \end{subequations}
    }
\end{model}

\subsection{Synchronous Generator Model}
To model the dynamic behavior of synchronous generators with sufficient fidelity, we adopt the \emph{classical machine model}~\cite{sauer1998power}, a simplified representation derived from the two-axis model. For each generator $g \in \mathcal{G}$, the full set of differential equations is given by:
\begin{align}
\small
\label{eq:2}
    \begin{bmatrix}
    \dot e^{\prime g}_q(t) \\
    \dot e^{\prime g}_d(t) \\
    \dot \delta^g(t) \\
    \dot \omega^g(t)
    \end{bmatrix}
    =
    \begin{bmatrix}
    0 \\
    0 \\
    \omega_s \left( \omega^g(t) - \omega_s \right) \\
    \frac{1}{m^g} \left( p_m^g - e^{\prime g}_d(t) i_d^g(t) - e^{\prime g}_q(t) i_q^g(t) - d^g \left( \omega^g(t) - \omega_s \right) \right)
    \end{bmatrix}
\end{align}
where $e^{\prime g}_d(t)$ and $e^{\prime g}_q(t)$ are the internal voltage components of the generator in the $D$-$Q$ reference frame, $\delta^g(t)$ is the rotor angle, and $\omega^g(t)$ is the rotor's angular frequency. The synchronous angular frequency is denoted by $\omega_s$, $p_m^g$ is the mechanical input power supplied by the turbine, $d^g$ is the damping coefficient, and $m^g$ is the inertia constant. $i_d^g(t)$ and $i_q^g(t)$ are the stator currents. They are computed in the machine's rotating reference frame as:
\begin{align}
\small
\label{eq:stator_current}
    \begin{bmatrix}
i_d^g(t) \\
i_q^g(t)
\end{bmatrix}
\! =\!
\begin{bmatrix}
0 \!\!\!&\!\!\! -x^{\prime g}_d \\
x^{\prime g}_d \!\!\!&\!\!\! 0
\end{bmatrix}^{-1}
\begin{bmatrix}
e^{\prime g}_d(t) - |V_g| \sin(\delta^g(t) - \theta_g) \\
e^{\prime g}_q(t) - |V_g| \cos(\delta^g(t) - \theta_g)
\end{bmatrix}
\end{align}
where, $x^{\prime g}_d$ is the transient reactance, $|V_g|$ is the voltage magnitude at the generator terminal, and $\theta_g$ is the corresponding phase angle. $x^{\prime g}_d$, the transient reactance, is a constant machine parameter. 

The state vector for generator $g \in \mathcal G$ is thus given by 
$\begin{bmatrix}
e^{\prime g}_q & e^{\prime g}_d & \delta^g & \omega^g
\end{bmatrix}^\intercal$. Note that in the classical machine model, $e^{\prime g}_q$ and $e^{\prime g}_d$ are held constant (i.e., $\dot e^{\prime g}_q = \dot e^{\prime g}_d = 0$), simplifying the system dynamics. Given this assumption, the generator model reduces to a second-order system:
\begin{subequations}
\small
    \begin{flalign}
        \dot \delta^g(t) &= \omega_s (\omega^g(t) - \omega_s) 
        \label{eq:delta}\\
        \dot \omega^g(t) &= \frac{1}{m^g} \left( p_m^g - d^g(\omega^g(t) - \omega_s) \right) - \frac{e^{\prime g}_{q}(0) |V_g|}{x^{\prime g}_d m^g} \sin(\delta^g(t) - \theta_g) 
        \label{eq:omega}
    \end{flalign}
    \label{eq:generator_dynamics}
\end{subequations}

We refer readers to~\cite{sauer1998power} for a comprehensive discussion of the classical model and its derivation from higher-order generator dynamics.

\subsection{Initial Values for the Generator States}
For each generator $g \in \mathcal{G}$, the initial values of the rotor angle $\delta^g(0)$ and the internal electromotive force (EMF) $e^{\prime g}_q(0)$ are computed based on the steady-state of the generator dynamics~\eqref{eq:generator_dynamics}, i.e., $\dot{\delta}^g(t) = 0$ and $\dot{\omega}^g(t) = 0$. Under this condition, the steady-state power balance equations for active and reactive power~\cite{machowski2020power} yield:
\begin{align}
\small
    &\frac{e^{\prime g}_q(0)\, |V_g| \sin(\delta^g(0) - \theta_g)}{x^{\prime g}_d} - p^r_g = 0, \label{eq:init_delta} \\
    &\frac{e^{\prime g}_q(0)\, |V_g| \cos(\delta^g(0) - \theta_g) - |V_g|^2}{x^{\prime g}_d} - q^r_g = 0. \label{eq:init_emf}
\end{align}

\noindent
These equations are solved simultaneously to determine the initial rotor angle $\delta^g(0)$ and EMF $e^{\prime g}_q(0)$ consistent with the given steady-state active and reactive power set-points $p^r_g$ and $q^r_g$, and terminal voltage magnitude and angle $|V_g|$ and $\theta_g$.

Following the same steady-state assumption, the initial rotor speed is given by:
\begin{align}
    \omega^g(0) = \omega_s. \label{eq:init_omega}
\end{align}

\subsection{Modeling Stability of Synchronous Generator Dynamics}
To ensure the dynamic stability of a synchronous generator $g \in \mathcal{G}$, it is essential to maintain its rotor angle $\delta^g(t)$ within a permissible range at all times. This condition reflects the requirement for \emph{rotor angle stability}, which refers to the ability of a generator to remain in synchronism with the rest of the system following small disturbances. Physically, rotor angle stability ensures that the generator's electromagnetic torque is balanced with the mechanical torque, preventing sustained angular divergence. 

A widely accepted criterion in the power systems literature~\cite{sauer1998power} is to impose an upper bound on the rotor angle to prevent loss of synchronism. This is typically enforced through the following constraint:
\begin{align}
    \delta^g(t) \leqslant \delta^{\max} \quad \forall t \geqslant 0, \label{eq:stab-limit}
\end{align}
where $\delta^{\max}$ is a system-specific stability threshold. This threshold is typically set based on transient stability margins, system protection settings, or empirical assessments of safe rotor-angle deviations under disturbances.

Violation of this constraint indicates potential instability, characterized by the generator's inability to resynchronize with the system, which may result in generator tripping or cascading failures. As such, enforcing rotor angle bounds is a fundamental requirement in stability-constrained dispatch and forms the basis for integrating dynamic constraints into the ACOPF framework.

\subsection{Stability-Constrained ACOPF}
\label{ssec:stability_constrained_OPF}
\begin{model}[tb]
    \caption{The Stability-Constrained ACOPF Problem}
    \label{model:stability_constrained_ac_opf}
    {\small
    \begin{subequations}
    \begin{align}
        \mbox{\bf variables:} \;\;
        & S^r_i, V_i \;\; \forall i\in \mathcal{N}, \;\; \delta^g(t), \omega^g(t) \;\; \forall g\in \mathcal{G}, \;\; \nonumber
         \\ & S_{ij}     \;\; \forall(i,j)\in \mathcal{L} \nonumber \\
        \mbox{\bf minimize:} \;\;
         & \sum_{i \in \mathcal{G}}  c_{2i} (\Re(S^r_i))^2 + c_{1i}\Re(S^r_i) +  c_{0i} \label{scac_obj} \\
        \textbf{subject to:} \;\;
        & \text{\eqref{eq:ac_1} -- \eqref{eq:ac_7} and} \nonumber \\ 
        & \text{\eqref{eq:generator_dynamics} -- \eqref{eq:stab-limit}} \;\; \forall g \in \mathcal G \nonumber 
    \end{align}
    \end{subequations}
    }
\end{model}
\smallskip
\noindent

The mathematical formulation for the stability-constrained ACOPF (SC-ACOPF) problem is shown in Model~\ref{model:stability_constrained_ac_opf}. This problem determines the optimal power dispatch subject to physical laws, operational constraints, and dynamic stability requirements, as captured by constraints~\eqref{eq:ac_1} -- \eqref{eq:ac_7} and \eqref{eq:generator_dynamics} -- \eqref{eq:stab-limit}.
Given the loads $\bm{S}^d$, the objective is to determine the dispatch and voltage variables $\bm{S}^r$ and $\bm{V}$ that minimize the cost function~\eqref{scac_obj} while satisfying both steady-state AC power flow equations and generator stability constraints.
The incorporation of synchronous generator dynamics introduces a significant modeling challenge due to the coupling between the OPF decision variables--namely, the complex power injection $S^r_g$ and terminal voltage $V_g$--and the generator state variables $\delta^g(t)$ and $\omega^g(t)$, as represented in constraints~\eqref{eq:generator_dynamics} -- \eqref{eq:stab-limit}. The nonlinear and time-dependent nature of the generator model further increases the dimensionality and complexity of the optimization problem.

Traditional discretization-based methods, which simulate the dynamic equations over time, become computationally prohibitive when integrated directly into the ACOPF framework. To address this, we develop a data-driven approach using Gaussian Process (GP) regression to approximate generator dynamics and their stability constraints in \eqref{eq:stab-limit}. This surrogate modeling strategy enables the efficient enforcement of stability constraints via chance constraints, thereby allowing the steady-state ACOPF formulation to account for dynamic behavior without incurring the cost of explicit time-domain simulation. The result is a computationally tractable framework that maintains both economic efficiency and compliance with stability.

\section{Gaussian Processes (GPs) - Preliminaries} 
\label{sec:gp-intro}
In this section, we review the preliminaries of using Gaussian Processes (GPs) for supervised regression tasks. Suppose we are interested in learning an unknown function $f: \mathbb{R}^n \to \mathbb{R}$ from a dataset of input-output pairs. We assume the observations follow the model:
$$
y_i = f(\bm x_i) + \epsilon_i, \quad \epsilon_i \sim \mathcal{N}(0, \sigma_n^2),
$$
where $\epsilon_i$ represents zero-mean Gaussian observation noise with variance $\sigma_n^2$, and $\bm x_i \in \mathbb{R}^n$ is the $i$-th input. A Gaussian Process
is a non-parametric prior over functions, fully specified by a mean function $m(\bm x)$ and a covariance function $k(\bm x, \bm x^\prime)$. Formally, we write:
\begin{gather*}
f(\bm x) \sim \mathcal{GP}(m(\bm x), k(\bm x, \bm x^\prime)), ~~\text{where} \\
m(\bm x) = \mathbb{E}_{f \sim \mathcal{GP}}[f(\bm x)], \text{ and } \\
k(\bm x, \bm x^\prime) = \mathbb{E}_{f \sim \mathcal{GP}}[(f(\bm x) - m(\bm x))(f(\bm x^\prime) - m(\bm x^\prime))].
\end{gather*}

Given a dataset $\mathcal{D} = \{(\bm x_i, y_i)\}_{i=1}^N$ with input matrix $\bm X = [\bm x_1, \dots, \bm x_N]^\intercal$ and output vector $\bm y = [y_1, \dots, y_N]^\intercal$, we define the joint prior distribution over the training outputs $\bm y$ and the function values $f_*$ at a set of test inputs $\bm X_*$ as:
\begin{equation}
\small
\begin{bmatrix}
\bm y \\
f_*
\end{bmatrix}
\sim \mathcal{N} \left(
\bm{0},
\begin{bmatrix}
K(\bm X, \bm X) + \sigma_n^2 \bm I & K(\bm X, \bm X_*) \\
K(\bm X_*, \bm X) & K(\bm X_*, \bm X_*)
\end{bmatrix}
\right),
\end{equation}
where, $K(\bm X, \bm X)$ denotes the kernel matrix evaluated at all pairs of training inputs, $K(\bm X_*, \bm X)$ the cross-covariance between test and training inputs, and $K(\bm X_*, \bm X_*)$ the covariance matrix of the test inputs.

The hyperparameters of the GP model--including kernel parameters such as length scale and signal variance--are typically optimized by maximizing the marginal log-likelihood of the observed data:
\begin{align}
\log p(\bm y \mid \bm X) =\,
& - \frac{1}{2} \bm y^\top (\bm K + \sigma_n^2 \bm I)^{-1} \bm y - \frac{1}{2} \log |\bm K + \sigma_n^2 \bm I|
  - \frac{N}{2} \log 2\pi
\end{align}
where, $\bm K \triangleq K(\bm X, \bm X)$.

At inference time, for a new test input $\bm x_* \in \mathbb{R}^n$, the predictive distribution for the function value $f_*$ is Gaussian with mean and variance given by:
\begin{align}
\mathbb{E}[f_* \mid \bm X, \bm y, \bm x_*] &= 
K(\bm x_*, \bm X) (\bm K + \sigma_n^2 \bm I)^{-1} \bm y, \\
\text{Var}[f_* \mid \bm X, \bm y, \bm x_*] &=
K(\bm x_*, \bm x_*) - K(\bm x_*, \bm X) (K + \sigma_n^2 \bm I)^{-1} K(\bm X, \bm x_*).
\end{align}
This closed-form Bayesian posterior yields point predictions and an explicit measure of model uncertainty. This is crucial in applications requiring robustness or risk-sensitive decision-making, particularly when model training is limited by data. We will use the uncertainty estimate in the next section to reformulate the GP-based surrogate stability models as chance constraints.
\begin{figure}[htbp]
    \centering
    \begin{subfigure}[b]{0.45\textwidth}
        \centering
        \includegraphics[scale=0.6]{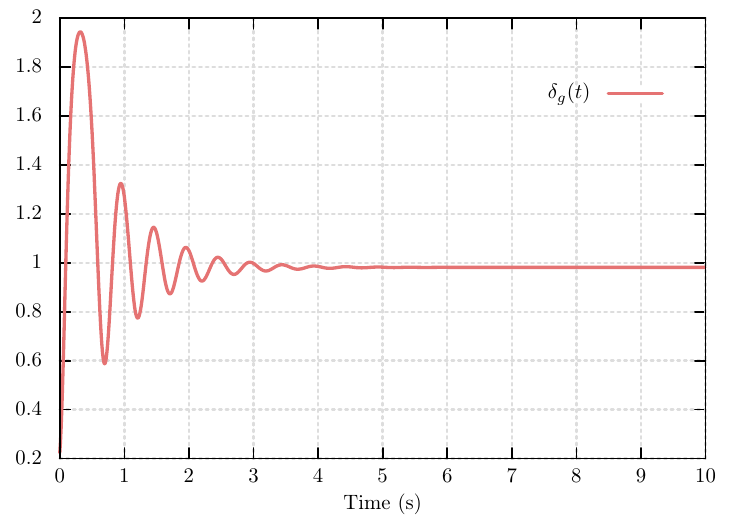}  
        \caption{Rotor angle when dynamics in \eqref{eq:generator_dynamics} is stable.}
        \label{fig:stable_traj}
    \end{subfigure}
    \hfill
    \begin{subfigure}[b]{0.45\textwidth}
        \centering
        \includegraphics[scale=0.6]{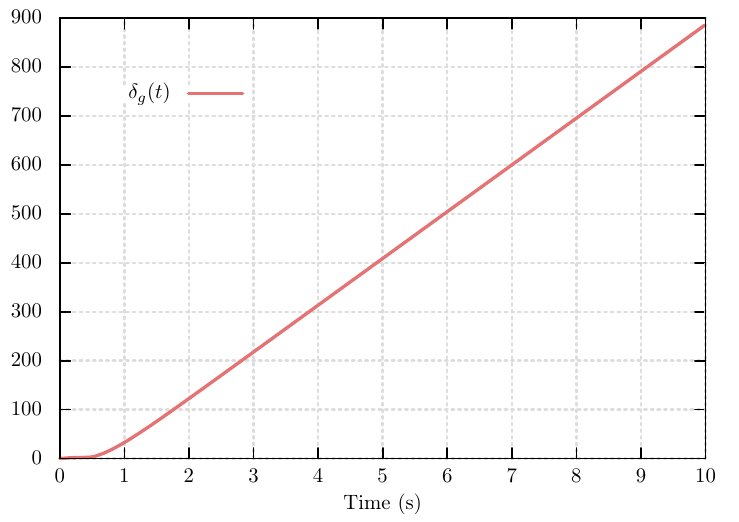}  
        \caption{Rotor angle when dynamics in \eqref{eq:generator_dynamics} is unstable.}
        \label{fig:unstable_traj}
    \end{subfigure}
    \caption{In stable operating conditions, the rotor angle $\delta^g(t)$ exhibits sinusoidal oscillations modulated by exponential decay. Under unstable conditions, $\delta^g(t)$ grows monotonically over time, typically exhibiting exponential or polynomial divergence.}
\end{figure}

\section{Design of GP-Based Stability Surrogates}
\label{sec:method}
\subsection{Qualitative View of Generator Dynamics} \label{subsec:qualitiative}
The dynamic behavior of synchronous generators, as described by the nonlinear differential equations in~\eqref{eq:generator_dynamics}, exhibits patterns that closely resemble those of second-order linear systems. Specifically, despite the inherent nonlinearity, the rotor angle trajectories qualitatively fall into two distinct categories: damped oscillations in stable regimes and unbounded growth in unstable regimes. This behavior can be observed in Figures~\ref{fig:stable_traj} and~\ref{fig:unstable_traj}, which show the system response under stable and unstable conditions, respectively, using initial conditions given by~\eqref{eq:init_delta} -- \eqref{eq:init_omega}.

In stable operating conditions, the rotor angle $\delta^g(t)$ exhibits sinusoidal oscillations modulated by exponential decay, reminiscent of under-damped second-order linear systems. In contrast, under unstable conditions, $\delta^g(t)$ grows monotonically over time, typically exhibiting exponential or polynomial divergence. This qualitative distinction is robust across a range of generator parameters characterizing the generator response in~\eqref{eq:generator_dynamics} such as mechanical input power $p_m^g$ and damping coefficient $d^g$, and operating points defined by $p^r_g$, $q^r_g$, $|V^g|$, and $\theta^g$. These parameters not only determine whether the generator operates in a stable or unstable regime but also influence the oscillation frequency and growth/decay rate of the rotor angle trajectory. Beyond this qualitative interpretation, we empirically analyzed rotor-angle trajectories obtained from extensive simulations across a wide range of operating conditions. Across these experiments, the trajectories consistently exhibit a dominant envelope behavior, either exponential-like decay in stable regimes or sustained growth in unstable conditions, modulated by an oscillatory component. This recurrent empirical pattern supports the use of a simplified surrogate that captures the overall growth or decay trend of the trajectory in the subsequent Section \ref{subsec:exp-design}.

From a control-theoretic perspective, this dichotomy in behavior aligns with classical Lyapunov-based notions of stability \cite{khalil2002nonlinear}. A stable generator trajectory remains confined within a bounded region of the state space, while an unstable one exhibits diverging trajectories, violating the boundedness condition. In SC-ACOPF formulations, this is encoded through the constraint~\eqref{eq:stab-limit}, which requires the rotor angle to remain below a specified threshold $\delta^{\max}$ for all $t \geqslant 0$. The clear separation in dynamic response types, along with their dependence on physical parameters and initial conditions, supports the use of simplified qualitative criteria for assessing stability within power system optimization.

\subsection{Exponential Surrogate Design} \label{subsec:exp-design}
Building on the qualitative behavior of generator dynamics during stable and unstable operations discussed in the previous section, we now introduce a compact surrogate model for approximating the trajectory of the rotor angle $\delta^g(t)$, as follows: 
\begin{align}
    \delta^g(t) \approx \alpha \cdot \exp(\beta t), \label{eq:exp_surrogate}
\end{align}
where $\alpha, \beta \in \mathbb{R}$ are case-specific parameters that depend on the operating point and system characteristics, and can be estimated using a nonlinear least squares fitting algorithm~\cite{markwardt2009nonlinearsquaresfittingidl}. 
We remark that the surrogate in \eqref{eq:exp_surrogate} is sufficient but not necessary, i.e., other approximations could, in principle, be used to capture the overall trend of the trajectory, but the exponential form offers a particularly direct and interpretable stability indicator. The motivation for this choice stems from the fact that, despite the dynamic model of a synchronous generator in \eqref{eq:generator_dynamics} governed by nonlinear second-order differential equations, the behavior around an operating point is well approximated by a linearized second-order system whose solutions take the form $\delta_g(t) = \exp{(\sigma t)} \sin (\omega t +\phi)$, where the exponential term, $\exp{(\sigma t)}$, governs the envelope of the response and dictates whether oscillations, given by $\sin(\omega t + \phi)$, decay or diverge. Since the goal of the surrogate is not to reconstruct trajectories pointwise but rather to extract a stability-relevant descriptor suitable for integration within the ACOPF formulation, capturing this exponential envelope is sufficient for accurate stability classification, as confirmed by the experimental results in Section \ref{sec:exp}. 
\begin{figure}[t!]
    \centering
    \begin{subfigure}[b]{0.45\textwidth}
        \centering
        \includegraphics[scale=0.6]{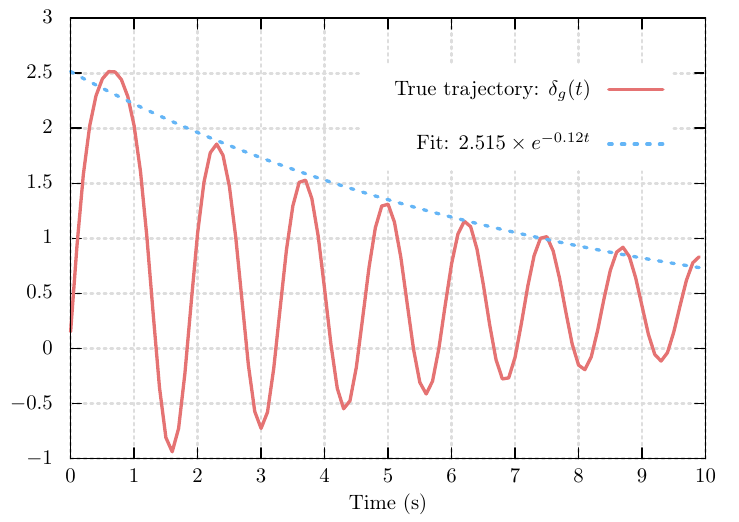}  
        \caption{$\delta_g(t)$ and its fit using \eqref{eq:exp_surrogate} for a stable dynamics.}
        \label{fig:stable_traj_with_envelope}
    \end{subfigure}
    \hfill
    \begin{subfigure}[b]{0.45\textwidth}
        \centering
        \includegraphics[scale=0.6]{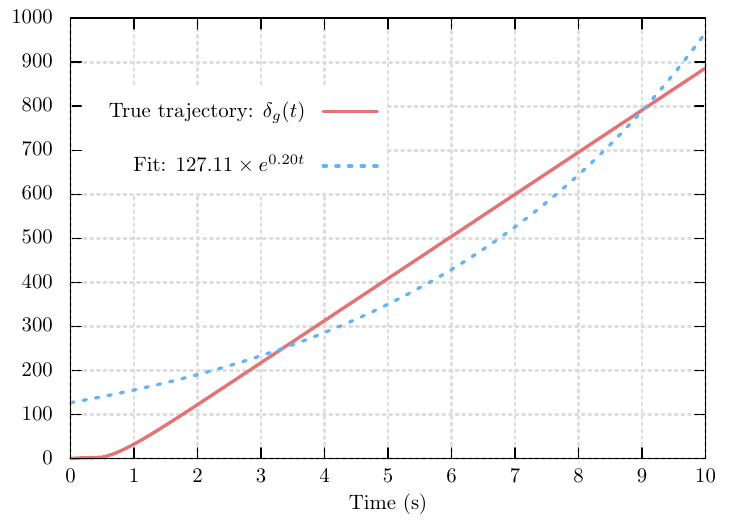}  
        \caption{$\delta_g(t)$ and its fit using \eqref{eq:exp_surrogate} for an unstable dynamics.}
        \label{fig:unstable_traj_with_envelope}
    \end{subfigure}
    \caption{The exponential surrogate fit effectively captures the growth or decay trend of the trajectory of $\delta_g(t)$. }
\end{figure}

Figures~\ref{fig:stable_traj_with_envelope} and~\ref{fig:unstable_traj_with_envelope} show rotor angle trajectories in stable and unstable regimes, discussed in Section~\ref{subsec:qualitiative}, but overlaid with the exponential surrogate~\eqref{eq:exp_surrogate}. The $\alpha$ and $\beta$ for the surrogate function are optimized using the nonlinear least squares routine to provide a close approximation of the envelope of the true dynamic behavior. The exponential surrogate in~\eqref{eq:exp_surrogate} is appealing for several reasons. First, it compactly captures the essential growth or decay trend of the rotor angle trajectory. 
The sign and magnitude of $\beta$ provide a direct interpretation: $\beta < 0$ corresponds to decaying (stable) trajectories, while $\beta > 0$ captures growing (unstable) dynamics. The magnitude indicates how large (small) the growth (decay) rate is in the dynamics. The coefficient $\alpha$ sets the initial response scale and can be tied to the system initialization or normalized if necessary. 

Second, from an optimization perspective, the surrogate form in~\eqref{eq:exp_surrogate} provides a tractable, differentiable representation. The simplicity of this functional form allows for efficient integration into probabilistic or chance-constrained frameworks within SC-ACOPF, which are required when generator stability must be enforced under uncertain inputs or model uncertainty.

\textcolor{black}{Finally, we note that the stability constraint in \eqref{eq:stab-limit} is an operational requirement; in our experiments, trajectories satisfying this condition also yielded $\beta < 0$. Under the exponential envelope approximation, both $\beta < 0$ and $\delta_g(t)\leqslant \delta^{\max}$ act as alternative constraints capturing the same underlying stability behavior. Accordingly, the $\beta$-based condition is best interpreted as a validated surrogate for stability-oriented dispatch, not as a formal replacement for all transient-security definitions.}

\subsection{GP Regression for Stability Prediction} \label{subsec:GP-regression}
We train a Gaussian Process (GP) regression model to estimate the exponential rate parameter $\beta$ that characterizes the stability of rotor angle trajectory $\delta^g(t)$ for each generator $g \in \mathcal{G}$. As introduced in Section~\ref{subsec:exp-design}, the generator dynamics exhibit a dominant exponential trend over time, which can be compactly captured using the surrogate form defined in~\eqref{eq:exp_surrogate}. 

To learn the mapping between operating conditions and the corresponding $\beta$, we construct a dataset $\mathcal{D} = \{ (\bm{x}_i, y_i) \}_{i=1}^{N}$. Each input vector $\bm{x}_i = (\delta^g(0), \omega^g(0), |\hat{V}_g|, \hat{\theta}_g)$ contains the initial rotor angle and angular speed, along with the terminal voltage magnitude and phase at the generator bus. 
The target value $y_i = \beta_i$ is computed by fitting the exponential surrogate to a simulated rotor angle trajectory $\delta^g(t)$, obtained by numerically solving the generator dynamics~\eqref{eq:generator_dynamics}. The fitting is performed using a nonlinear least squares routine~\cite{markwardt2009nonlinearsquaresfittingidl}.

To ensure that the GP model generalizes across the full space of feasible operating conditions, the training data are generated by sampling operating points within practical bounds. Specifically, the voltage magnitude $|\hat{V}_g|$ and phase angle $\hat{\theta}_g$ are sampled from uniform distributions over intervals defined by the operational constraints~\eqref{eq:ac_1} and~\eqref{eq:ac_6}, i.e., $|\hat{V}_g| \sim U(v^l_g, v^u_g)$ and $\hat{\theta}_g \sim U(\theta^l_g, \theta^u_g)$. Similarly, the active and reactive power injections $p^r_g$ and $q^r_g$ are sampled from their respective feasible regions defined by~\eqref{eq:ac_2}. Given these sampled values, the corresponding initial conditions $\delta^g(0)$ and $\omega^g(0)$ are computed using the initialization formulas~\eqref{eq:init_delta} -- \eqref{eq:init_omega}.

The GP regression model is then trained to learn the mapping:
\begin{align}
    \beta = f(\delta^g(0), \omega^g(0), |\hat{V}_g|, \hat{\theta}_g), \label{eq:gp_regression}
\end{align}
where $f$ is modeled as a Gaussian Process using the general technique presented in Section \ref{sec:gp-intro}. This formulation enables the GP to predict the exponential trend in rotor-angle dynamics under new operating conditions, along with an associated uncertainty measure, which is subsequently used to formulate chance constraints in the SC-ACOPF. Although the terminal voltage magnitude, $|\hat{V}_g|$, and phase,  $\hat{\theta}_g$,  are not dynamic state variables in the generator model, they directly influence both the initialization and subsequent evolution of the rotor angle dynamics. Importantly, different ACOPF operating points correspond to different voltage magnitudes, phase angles, and power injections, which in turn induce distinct dynamic regimes. Hence, to enable the GP surrogate to distinguish among these operating conditions, it is necessary to include $|\hat{V}_g|$ and $\hat{\theta}_g$ explicitly as input features. Without these variables, the GP would be unable to differentiate operating points that share identical initial states $\delta^g(0)$ and $\omega^g(0)$ but evolve differently due to distinct bus voltages. Including these features ensures that the learned mapping in~\eqref{eq:gp_regression} captures the dependence of stability on the operating conditions determined by the ACOPF solution.

This training framework allows each GP model to approximate a family of generator behaviors across a wide range of feasible operating points. Once trained, the GP models serve as efficient surrogates for assessing rotor angle stability and can be seamlessly integrated into the SC-ACOPF formulation. The next section discusses how this integration is achieved using chance constraints.

\subsection{Integration into ACOPF via Chance Constraints} \label{subsec:cc}
To incorporate the GP-based surrogate stability model into the ACOPF framework, we impose the condition that the exponential rate parameter $\beta$--as introduced in the surrogate model~\eqref{eq:exp_surrogate}--remains non-positive to ensure bounded rotor angle behavior over time. This requirement is formalized as:
\begin{align}
    \beta \leqslant 0. \label{eq:stab-limit-surrogate}
\end{align}
Since $\beta$ is predicted using a Gaussian Process regression model trained on generator operating conditions, the GP provides not just a point estimate, but a full predictive distribution for $\beta$, capturing epistemic uncertainty due to limited training data and surrogate approximation error. Specifically, for a given generator $g \in \mathcal{G}$ and operating point $\bm{x}_g = (\delta^g(0), \omega^g(0), |\hat{V}_g|, \hat{\theta}_g)$, the GP yields a predictive mean $\mu_g(\bm{x}_g)$ and standard deviation $\sigma_g(\bm{x}_g)$, both of which are smooth, nonlinear functions of the input features and the learned model parameters.
To ensure stability with high probability while accounting for model uncertainty, we introduce a chance constraint:
\begin{align}
    \mathbb{P}(\beta \leqslant 0) \geqslant 1 - \varepsilon, \label{eq:prob-stab-limit}
\end{align}
where $\varepsilon \in (0,1)$ is a user-defined risk tolerance. This constraint ensures that the generator remains in a stable regime with probability at least $1 - \varepsilon$ under the GP's predictive distribution.
Because the GP prediction for $\beta$ follows a Gaussian distribution, $\beta \sim \mathcal{N}(\mu_g(\bm{x}_g), \sigma_g^2(\bm{x}_g))$, the chance constraint~\eqref{eq:prob-stab-limit} can be equivalently reformulated as the deterministic inequality:
\begin{align}
    \mu_g(\bm{x}_g) + \Phi^{-1}(1 - \varepsilon) \cdot \sigma_g(\bm{x}_g) \leqslant 0, \label{eq:deterministic-chance}
\end{align}
where $\Phi^{-1}$ is the inverse cumulative distribution function of the standard normal distribution. This constraint is linear in $\mu_g(\bm{x}_g)$ and $\sigma_g(\bm{x}_g)$ and can be efficiently handled by standard nonlinear optimization solvers.

A key advantage of using GP models in this setting is their compatibility with gradient-based solvers such as IPOPT~\cite{Wchter2006OnTI}. The GP inference step involves only matrix-vector operations, and both $\mu_g(\bm{x}_g)$ and $\sigma_g(\bm{x}_g)$ are differentiable functions of the generator operating point. This eliminates the need to simulate differential equations at each iteration, as is required in discretization-based dynamic formulations, thereby significantly reducing the computational burden. Fig.~\ref{fig:flow-chart} presents a schematic overview of the proposed GP-SC-ACOPF framework, highlighting the offline training of the proposed GP surrogates and their direct embedding as chance constraints within the SC-ACOPF optimization problem.

\begin{figure}[!ht] 
\centering 
\begin{tikzpicture}[
    node distance=2.5cm,
    every node/.style={font=\small},
    block/.style={
        rectangle, draw, fill=green!10,
        text width=15em, text centered,
        rounded corners, minimum height=3.2em
    },
    block2/.style={
        rectangle, draw, fill=blue!15,
        text width=15em, text centered,
        rounded corners, minimum height=3.2em
    },
    line/.style={draw, -{Latex}}
]


\node[align=center] (offline_title)
{\textbf{Offline Phase}\\
\footnotesize For each $g\in\mathcal{G}$ train a GP};

\node[block, below=0.6cm of offline_title] (sample)
{Sample operating points\\
$|V_g|,\theta_g,p_g^r,q_g^r$};

\node[block, below=.5cm of sample] (simulate)
{Compute initial conditions \eqref{eq:init_delta}--\eqref{eq:init_omega}
and solve generator dynamics in Eq. \eqref{eq:generator_dynamics}};

\node[block, below=.5cm of simulate] (fit)
{Fit $\delta^g(t)\approx \alpha e^{\beta t}$ and
extract $\beta$};

\node[block, below=.5cm of fit] (dataset)
{Construct dataset $\mathcal{D}=\{(\bm x_i, y_i)\}$\\
$\bm x_i=(\delta^g(0),\omega^g(0),|V_g|,\theta_g)$, $y_i = \beta_i$};

\node[block, below=.5cm of dataset] (train)
{Train a GP for $\beta=f(\delta^g(0),\omega^g(0),|V_g|,\theta_g)$\\
Obtain predictive $\mu_g(\bm x_g),\sigma_g(\bm x_g)$};


\node[align=center, right=3.2cm of offline_title] (online_title)
{\textbf{Online Phase}\\
\footnotesize Stability-constrained dispatch};

\node[block2, below=0.5cm of online_title] (formulation)
{Formulate GP-SC-ACOPF:\\
\small
- AC power flow (1b)--(1h)\\
- GP-based stability constraints};

\node[block2, below=.5cm of formulation] (gpconstraint)
{\small Embedded GP stability constraint:\\
For each $g\in\mathcal{G}$, enforce\\
$\mu_g(\bm x_g) + \Phi^{-1}(1-\varepsilon)\,\sigma_g(\bm x_g) \leqslant 0$};

\node[block2, below=.5cm of gpconstraint] (solve)
{Solve with Ipopt the \\ GP-SC-ACOPF problem};

\node[block2, below=.5cm of solve] (dispatch)
{Stability-aware optimal dispatch};


\path[line] (sample) -- (simulate);
\path[line] (simulate) -- (fit);
\path[line] (fit) -- (dataset);
\path[line] (dataset) -- (train);

\path[line] (formulation) -- (gpconstraint);
\path[line] (gpconstraint) -- (solve);
\path[line] (solve) -- (dispatch);

\path[line] (train.east) -- ++(0.8,0) |- (formulation.west);

\end{tikzpicture}
\caption{Methodology of the proposed GP-SC-ACOPF framework. The offline phase trains one GP per generator, and the online phase integrates the surrogates into ACOPF via chance constraints.} 
\label{fig:flow-chart} 
\end{figure}

\section{Experiments}
\label{sec:exp}
Having established the framework integrating GP-based generator dynamics within ACOPF, this section focuses on assessing it on three key aspects:
\begin{enumerate}[leftmargin=*,nosep]
    \item The predictive accuracy of the proposed GP models for synchronous-generator dynamics in realistic power-grid scenarios, benchmarked against alternative learning-based surrogates for differential equations.
    \item The operational efficacy of GP-SC-ACOPF, the GP-based stability-constrained AC optimal power flow, formulation relative to the conventional steady-state ACOPF in ensuring dynamic stability. We report statistics on the proportion of cases in which each method satisfies stability requirements under varying load conditions and generator-contingency scenarios, across three representative test systems.
    \item The sensitivity of GP-SC-ACOPF dispatch outcomes to the choice of confidence level~$\varepsilon$ in the chance constraints \eqref{eq:deterministic-chance}. We analyze how adjustments to~$\varepsilon$ affect the total generation cost and compare them with the steady-state ACOPF solutions.
\end{enumerate}

All the implementations were done using the Python programming language; the ACOPF models were constructed using the \texttt{Egret} package and solved with \texttt{IPOPT}, which is provided by the \texttt{Pyomo} package~\cite{bynum2021pyomo}. The GP models were trained using \texttt{GPYTorch}~\cite{gardner2018gpytorch} and integrated within the ACOPF model using \texttt{ROGP}~\cite{Wiebe2020robust}. 
The other learning surrogates used for benchmarking in the forthcoming section were implemented using \texttt{PyTorch}~\cite{paszke2019pytorchimperativestylehighperformance} and the \texttt{torchdiffeq} library. The solutions to the dynamic generator models were computed using the \texttt{SciPy} package. All experiments were carried out on a MacBook Pro equipped with an Apple M2 Pro chip and 16GB of RAM.

\subsection{Predictive Performance of the Gaussian Process Models}
\label{subsec:predictive-performance}
We first assess the ability of our Gaussian process (GP) surrogates to capture synchronous-generator dynamic stability via the exponential rotor-angle surrogate introduced in Section~\ref{subsec:exp-design}.  Each generator--though governed by the common model in~\eqref{eq:generator_dynamics}--exhibits unique inertia and damping parameters; our evaluation spans multiple generator parameterizations in standard IEEE test systems (IEEE 39-bus, 57-bus, and 118-bus). The choice of these three systems is purely motivated by the availability of open-source data for generator parameters. 
The procedure for generating the dataset $\mathcal D$ for this study is provided in Section \ref{subsec:GP-regression}. 

For benchmarking, we compare against three learning-based surrogates:
\begin{itemize}[leftmargin=*,nosep]
  \item Feed-forward neural networks (FNNs)~\cite{lecun2015deep} trained to predict $\beta$ directly from the $\mathcal D$ by minimizing mean squared error.
  \item Neural ODEs~\cite{chen2019neural} that learn a parametric vector field and generate $\hat\delta^g(t)$ by numerical integration directly from the input vector $\bm x_i$; we then fit the exponential surrogate to $\hat\delta^g(t)$ to extract $\hat\beta$.
  \item Physics-informed neural networks (PINNs)~\cite{RAISSI2019686} that ingest $\bm x_i$ and provides a pointwise $\hat\delta^g(t_k)$ over a uniform grid $\{t_k\}_k$, from which $\hat\beta$ is obtained via the same nonlinear least-squares routine.
\end{itemize}

All methods are evaluated as binary classifiers on the sign of $\beta$: a prediction is deemed correct if $\operatorname{sign}(\hat\beta)=\operatorname{sign}(\beta)$.  To mimic practical data scarcity, we train each model on datasets of $100$ and $1000$ samples and report accuracy on a fixed held-out test set of $1000$ examples.  This low-data regime highlights the sample efficiency and generalization capacity of the GP surrogates relative to the neural baselines.

\begin{table*}[htbp]
\small
    \centering
    \caption{Predicting the generator dynamics: accuracy (out of $1000$) of different methods across test cases and training sample sizes in the low sample regime, when interpreting the regression task as a binary classification task.}
    \label{tab:comparison_results}
    \begin{tabular}{lccccc}
        \toprule
        Test Case & Samples & PINN & Neural ODE & FNN & GP \\
        \midrule
        \multirow{2}{*}{\textbf{39-bus}}  & 100  & $0.621 \pm 0.085$  & $0.493 \pm 0.059$  & $0.532 \pm 0.116$  & $\bm{0.810 \pm 0.065}$ \\
        & 1000 & $0.571 \pm 0.063$  & $0.625 \pm 0.084$  & $0.544 \pm 0.043$  & $\bm{0.882 \pm 0.026}$ \\
        \midrule
        \multirow{2}{*}{\textbf{57-bus}}  & 100  & $0.604 \pm 0.102$  & $0.521 \pm 0.067$  & $0.662 \pm 0.132$  & $\bm{0.896 \pm 0.069}$ \\
        & 1000 & $0.567 \pm 0.067$  & $0.599 \pm 0.072$  & $0.693 \pm 0.038$  & $\bm{0.953 \pm 0.019}$ \\
        \midrule
        \multirow{2}{*}{\textbf{118-bus}} & 100  & $0.408 \pm 0.090$  & $0.357 \pm 0.085$  & $0.520 \pm 0.116$  & $\bm{0.872 \pm 0.071}$ \\
        & 1000 & $0.398 \pm 0.061$  & $0.560 \pm 0.053$  & $0.529 \pm 0.046$  & $\bm{0.893 \pm 0.033}$ \\
        \bottomrule
    \end{tabular}
\end{table*}

Table~\ref{tab:comparison_results} reports the mean and standard deviation of classification accuracy ($0$ to $1$, where $0.5$ indicates random guessing) for each surrogate across the IEEE 39-, 57-, and 118-bus systems. The GP surrogate consistently outperforms all baselines. In the 39-bus system with $100$ samples, GPs achieve $0.810 \pm 0.065$, outperforming the best baseline--PINN ($0.621 \pm 0.085$)--by $18.9\%$; with $1000$ samples, GP accuracy rises to $0.882 \pm 0.026$, a $25.7$-point advantage over PINN ($0.625 \pm 0.084$). Similarly, in the 57-bus network, GPs attain $0.896 \pm 0.069$ (100 samples) and $0.953 \pm 0.019$ (1000 samples), exceeding the top neural baseline (FNN at $0.662 \pm 0.132$ and $0.693 \pm 0.038$, respectively) by $23.4$ and $26.0$ points. In the 118-bus case, GPs achieve $0.872 \pm 0.071$ and $0.893 \pm 0.033$, outpacing their best competitor by $35.2$ and $33.3$ points under the two regimes.

These results in Table \ref{tab:comparison_results} demonstrate the superior sample efficiency and robustness of GP surrogates for dynamic-stability prediction, particularly in low-data settings, and underscore their suitability for integration within ACOPF formulations when simulation data is scarce.

\subsection{Stability under Varying Load Conditions}
\label{subsec:stability-loads}
This section demonstrates the benefit of solving the GP-SC-ACOPF, in lieu of a conventional steady-state ACOPF without stability constraints, to achieve a higher percentage of dynamically stable operating points across a wide range of load conditions.  We conduct experiments on three standard IEEE test cases (39-, 57-, and 118-bus systems), each featuring different numbers of generators and distinct inertia/damping parameters. For each system, individual bus loads are independently varied between 50\% and 150\% of their nominal values. Any randomized load profile for which the ACOPF solver terminates with any status other than local optimal is discarded, and another load profile is regenerated until feasibility is achieved.  
For this study, the value of $\varepsilon$ is set to $0.05$. 

For each test case, we collect $1000$ load realizations.  For each realization, we solve both the standard ACOPF and the GP-SC-ACOPF. We then evaluate dynamic stability by numerically integrating the generator dynamics in~\eqref{eq:generator_dynamics} with SciPy's ODE solver.  If all rotor-angle trajectories remain bounded (i.e., do not diverge), the dispatch is labeled \emph{stable}; otherwise, it is labeled \emph{unstable}.  Tables~\ref{tab:IEEE39_succed_fail}--\ref{tab:IEEE118_succed_fail} summarize the joint outcomes.

\noindent In the 39-bus case (Table \ref{tab:IEEE39_succed_fail}), the steady-state ACOPF yields a stable operating point in only $45.1\%$ of scenarios (sum of its ``stable-unstable'' and ``stable-stable'' rows), whereas the GP-SC-ACOPF achieves stability in $83.4\%$ of cases (sum of ``unstable-stable'' and ``stable-stable'' rows).  Notably, $42.9\%$ of cases that would diverge under ACOPF become stable under GP-SC-ACOPF, while only $4.6\%$ of ACOPF‐stable points lose stability when using the GP‐constrained formulation.

\begin{table}[htbp]
\small
\centering
\begin{minipage}{.3\textwidth}
\begin{tabular}{llr}
\toprule
ACOPF & \shortstack{GP-SC\\-ACOPF} & Percentage \\
\midrule
unstable & unstable & $12.0\%$ \\
\rowcolor{green!15} unstable & stable   & $42.9\%$ \\
\rowcolor{blue!10} stable   & unstable & $4.6\%$  \\
stable   & stable   & $40.5\%$ \\
\bottomrule
\end{tabular}
\caption{Stability for the 39-bus system.}
\label{tab:IEEE39_succed_fail}
\end{minipage}
\hfill
\begin{minipage}{.3\textwidth}
\begin{tabular}{llr}
\toprule
ACOPF & \shortstack{GP-SC\\-ACOPF} & Percentage \\
\midrule
unstable & unstable & $7.0\%$  \\
\rowcolor{green!15} unstable & stable   & $29.9\%$ \\
\rowcolor{blue!10} stable   & unstable & $2.9\%$  \\
stable   & stable   & $60.2\%$ \\
\bottomrule
\end{tabular}
\caption{Stability for the 57-bus system.}
\label{tab:IEEE57_succed_fail}
\end{minipage}
\hfill
\begin{minipage}{.3\textwidth}
\begin{tabular}{llr}
\toprule
ACOPF & \shortstack{GP-SC\\-ACOPF}& Percentage \\
\midrule
unstable & unstable & $11.0\%$ \\
\rowcolor{green!15} unstable & stable   & $30.7\%$ \\
\rowcolor{blue!10} stable   & unstable & $3.5\%$  \\
stable   & stable   & $54.8\%$ \\
\bottomrule
\end{tabular}
\caption{Stability for the 118-bus system.}
\label{tab:IEEE118_succed_fail}
\end{minipage}
\end{table}

\noindent For the 57-bus network, the base ACOPF is stable in $63.1\%$ of cases, while GP-SC-ACOPF attains stability in $90.1\%$ of scenarios.  The GP‐based approach rescues $29.9\%$ of cases that would otherwise diverge, with only a $2.9\%$ incidence of stability loss relative to the standard ACOPF solution.
\noindent In the largest 118-bus test case, ACOPF without stability constraints is stable only $58.3\%$ of the time, whereas GP-SC-ACOPF achieves $85.5\%$ stability.  GP‐SC-ACOPF converts $30.7\%$ of unstable ACOPF solutions into stable ones, at the cost of only $3.5\%$ of previously stable dispatch points. 

\textcolor{black}{In addition to reporting how often GP-SC-ACOPF stabilizes ACOPF-unstable operating points, it is equally important to quantify how often it produces unstable operating points that are instead stable under the standard ACOPF. To make this trade-off explicit, we reinterpret the joint stability results in Tables 2–4 in conditional form. For each benchmark, we evaluate (i) the probability that GP-SC-ACOPF recovers stability when ACOPF is unstable, and (ii) the probability that it introduces instability when ACOPF is stable. Across the three systems, the proposed method recovers stability in a large majority of ACOPF-unstable cases: 78.1\% (39-bus), 81.0\% (57-bus), and 73.6\% (118-bus). In contrast, transitions from stable to unstable operation remain limited, occurring in only 10.2\%, 4.6\%, and 6.0\% of ACOPF-stable cases, respectively. These results show that, although some degradation cases occur, they are significantly less frequent than the observed stability improvements. The net effect is therefore strongly positive, which is reflected in the overall stability rates increasing from 45.1\% to 83.4\% (39-bus), 63.1\% to 90.1\% (57-bus), and 58.3\% to 85.5\% (118-bus).}

Overall, these results confirm that incorporating GP‐based stability constraints dramatically increases the fraction of dynamically stable operating points across all loading conditions and network sizes.  The GP‐SC‐ACOPF yields robust stability improvements--often by more than $30\%$--while incurring only minimal cases of unnecessary conservatism.  

\begin{figure}[htbp]
    \centering
    \includegraphics[scale=0.7]{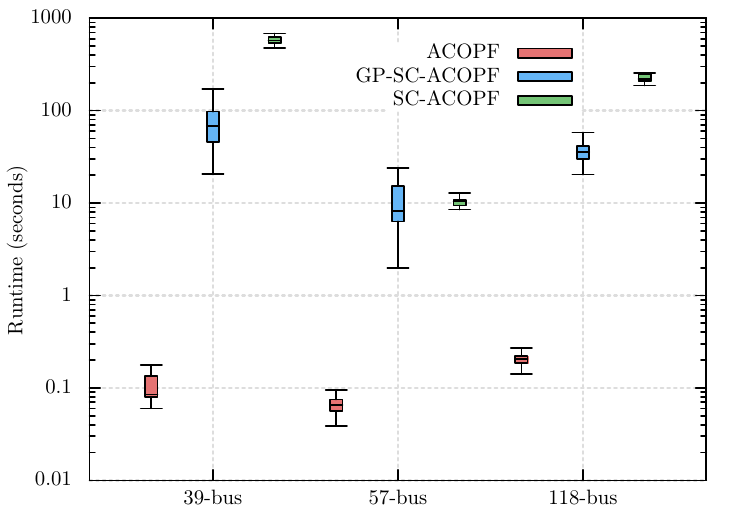}
    \caption{The time taken in seconds to solve the GP-SC-ACOPF, the SC-ACOPF and the ACOPF.}
    \label{fig:runt-time}
\end{figure}

Finally, it is important to note that improved stability does not come for free: as shown in Figure~\ref{fig:runt-time}, the average runtime of GP-SC-ACOPF is up to two orders of magnitude higher than that of the standard ACOPF. This trade-off between enhanced reliability and computational overhead must be carefully weighed in practical implementations.
Also shown in Figure~\ref{fig:runt-time} is the time required to solve the stability-constrained ACOPF in Model~\ref{model:stability_constrained_ac_opf} directly by discretizing the generator dynamics and solving the resulting discretized problem as a nonlinear optimization problem. Although this approach guarantees stability, the computation time for solving the SC-ACOPF can be up to three orders of magnitude higher than that of the standard ACOPF.
These results clearly demonstrate the substantial computational burden of the direct approach and underscore the need for data-driven methods that explicitly incorporate dynamic stability constraints while maintaining tractability.

\textcolor{black}{Given the negligible computational cost of solving ACOPF without stability constraints, a practical workflow is to first compute the ACOPF dispatch and invoke GP-SC-ACOPF only when the resulting operating point is dynamically unstable. Such two-stage approaches, in which model complexity is introduced selectively to improve feasibility, are standard in power system studies, including generator and load interconnection analyses \cite{GormanEtak2025_Joule}. The empirical results support this strategy. As shown in Tables~\ref{tab:IEEE39_succed_fail}--\ref{tab:IEEE118_succed_fail}, approximately 45--65\% of operating points are already dynamically stable under standard ACOPF, indicating that a substantial fraction of cases would not require the additional complexity of GP-SC-ACOPF. For the remaining unstable cases, GP-SC-ACOPF recovers stability in most scenarios, as evidenced by the observed transition statistics. While a small fraction of cases (on the order of 5--10\%) exhibit transitions from stable ACOPF solutions to unstable GP-SC-ACOPF outcomes, these are significantly outnumbered by unstable-to-stable corrections, leading to a net improvement in overall stability. Although this two-stage workflow is not explicitly benchmarked as an end-to-end pipeline, an evaluation that lies beyond the scope of this work, the reported results provide consistent empirical support for its effectiveness as a practical deployment strategy, balancing computational efficiency and dynamic security.}


\subsection{Stability under Generator Contingencies}
\label{subsec:contingency-stability}
This section evaluates the ability of GP-SC-ACOPF to deliver dynamically stable operating points under random N-1 generator contingencies. We adopt the same experimental setup as in Section~\ref{subsec:stability-loads}, but focus exclusively on the IEEE 118-bus system due to its larger number of generators. For each of the $1000$ feasible load realizations, we (i) solve the nominal dispatch via either the standard ACOPF or the GP-SC-ACOPF; (ii) randomly outage one generator; (iii) redistribute its lost generation proportionally among the remaining units; and (iv) solve a power flow to update bus voltages. We assess dynamic stability by integrating the generator model in~\eqref{eq:generator_dynamics} with SciPy's ODE solver. Similar to the results in Section \ref{subsec:stability-loads}, an outcome is labeled \emph{stable} if all generator rotor-angle trajectories remain bounded, and \emph{unstable} otherwise. Table~\ref{tab:IEEE118_gen_contingencies} summarizes the joint stability outcomes.

\begin{table}[htbp]
\small
\centering
\caption{Stability outcomes for the 118-bus system under random $N$-$1$ generator contingencies ($1000$ contingency scenarios).}
\label{tab:IEEE118_gen_contingencies}
\begin{tabular}{llr}
\toprule
ACOPF & GP-SC-ACOPF & Percentage \\
\midrule
unstable & unstable & $20.9\%$ \\
\rowcolor{green!15} unstable & stable   & $28.4\%$ \\
\rowcolor{blue!10} stable   & unstable & $12.0\%$ \\
stable   & stable   & $38.6\%$ \\
\bottomrule
\end{tabular}
\end{table}

In the ACOPF case, only $50.6\%$ of the contingency scenarios remain stable, whereas GP-SC-ACOPF achieves stability in $67.0\%$ of cases. Notably, GP-SC-ACOPF ``rescues'' $28.4\%$ of scenarios that would diverge under standard ACOPF, at the expense of a $12.0\%$ incidence where an ACOPF-stable dispatch becomes unstable under the GP-constrained formulation. These results highlight that the GP-SC-ACOPF significantly improves resilience to single‐generator outages. Note that using the mentioned workflow of ACOPF followed by GP-SC-ACOPF for deciding generation set points would provide stable solutions to $79.1\%$ of the generator contingency scenarios, which is significantly higher than the one with ACOPF alone.

\subsection{Stability under Line Contingencies}
\label{subsec:line-contingency-stability}
This section evaluates the ability of GP-SC-ACOPF to generate stable operating points under random $N$-$1$ line contingencies in the IEEE 118-bus system. The experimental setup mirrors that of Section~\ref{subsec:contingency-stability}, with the difference that contingencies are applied to transmission lines rather than generators. 
Table~\ref{tab:IEEE118_line_contingencies_rcd} summarizes the joint stability outcomes.
\begin{table}[htbp]
\small
\centering
\caption{Stability outcomes for the 118-bus system under $N$-$1$ line contingencies.}
\label{tab:IEEE118_line_contingencies_rcd}
\begin{tabular}{llr}
\toprule
ACOPF & GP-SC-ACOPF & Percentage \\
\midrule
unstable & unstable & $26.1\%$ \\
\rowcolor{green!15} unstable & stable   & $27.5\%$ \\
\rowcolor{blue!10} stable   & unstable & $10.5\%$ \\
stable   & stable   & $35.9\%$ \\
\bottomrule
\end{tabular}
\end{table}
Under standard ACOPF, $46.4\%$ of post-contingency operating points remain dynamically stable. In contrast, GP-SC-ACOPF achieves stability in $63.4\%$ of cases. Notably, GP-SC-ACOPF stabilizes $27.5\%$ of scenarios that would otherwise diverge under ACOPF, while introducing instability in $10.5\%$ of cases that were originally stable.
Overall, the net stabilization effect remains positive, with the GP-SC-ACOPF method delivering a substantial increase in resilience to transmission-line outages compared to the ACOPF baseline.
As in the generator-contingency setting, combining the two workflows (i.e., solving ACOPF and subsequently verifying or refining with GP-SC-ACOPF) yields stable solutions in $73.4\%$ of contingency scenarios, significantly exceeding the stability rate achieved by ACOPF alone.

\subsection{Sensitivity of Objective Value with respect to $\varepsilon$}
\label{subsec:objective_gap}
In addition to stability improvements, it is important to quantify the cost of enforcing the GP‐based chance constraints in terms of the objective value.  We define the relative objective gap as
\[
\text{Gap}(\%) \;=\; 100\%\;\times\;\frac{|J_{\rm GP\text{-}SC\text{-}ACOPF} - J_{\rm ACOPF}|}{J_{\rm ACOPF}},
\]
where \(J_{\rm GP\text{-}SC\text{-}ACOPF}\) and \(J_{\rm ACOPF}\) denote the optimal costs of the GP‐constrained and standard ACOPF, respectively.  Table~\ref{tab:objective_gap} consolidates these gaps for the 39-, 57-, and 118-bus systems under five values of $\varepsilon \in \{0.1, 0.05, 0.02, 0.01, 0.001\}$. 

\begin{table*}[htbp]
\small
    \centering
        \caption{Relative objective gaps (\%) of GP-SC-ACOPF compared to standard ACOPF for varying $\varepsilon$ values across IEEE test systems.}
    \label{tab:objective_gap}
    \begin{tabular}{lrrrrr}
        \toprule
        Test Case & $\varepsilon = 0.1$ & $\varepsilon = 0.05$ &$\varepsilon = 0.02$ & $\varepsilon = 0.01$ & $\varepsilon = 0.001$ \\
        \midrule
        39-bus  & $1.268$ & $1.274$ & $1.281$ & $1.281$ & \(\bm{1.290}\) \\
        57-bus  & $0.140$ & $0.169$ & $0.147$ & \(\bm{0.699}\) & $0.452$ \\
        118-bus & $0.715$ & $0.737$ & $0.745$ & $0.788$ & \(\bm{0.811}\) \\
        \bottomrule
    \end{tabular}
\end{table*}

As shown in Table~\ref{tab:objective_gap}, the additional cost of enforcing dynamic‐stability constraints remains modest across all networks.  In the 39-bus system, the gap increases only slightly from $1.268\%$ at $\varepsilon = 0.1$ to $1.290\%$ at $0.001$.  The 57-bus case exhibits a small peak of $0.699\%$ at $\varepsilon = 0.01$ confidence, while the 118-bus network sees a gradual rise from $0.715\%$ to $0.811\%$ over the same range.  These results indicate that the GP‐SC‐ACOPF can deliver substantial stability gains at a very limited increase in generation cost.  

\section{Conclusion and Future Work}
\label{sec:conclusion}
In this paper, we introduced a GP-based SC-ACOPF formulation that explicitly enforces dynamic‐stability constraints via chance constraints derived from the underlying GP estimates.  Our experimental evaluation on the IEEE 39-, 57-, and 118-bus systems demonstrates three key findings: (i) The GP surrogates achieve classification accuracies in excess of $87\%$ even in a low‐data regime ($100$ training samples), consistently outperforming PINNs, Neural ODEs, and FNNs by margins of up to $35\%$. (ii) Across thousands of randomized load scenarios, GP‐SC‐ACOPF is able to compute operating points that are stable in $30\%$ of cases that diverge under the standard ACOPF, raising the fraction of stable dispatches by more than $30\%$ in all test systems.  Under random $N$-$1$ generator outages, the proposed approach similarly converts nearly $30\%$ of potential failures into stable solutions. (iii) Enforcing GP‐based stability constraints incurs only a small objective‐value penalty (below $1.3\%$ in all networks) but increases solve times by roughly one order of magnitude.  This trade‐off between reliability and computational overhead must be carefully balanced from a practical standpoint. Together, these results establish that GP‐SC‐ACOPF provides a robust and sample‐efficient path to incorporating dynamic‐stability guarantees within the ACOPF framework, with negligible cost impact and manageable computational requirements.  

While our work lays a solid foundation, several avenues remain for further exploration: (i) Investigate decomposition techniques or advanced interior‐point heuristics to reduce the runtime gap between GP‐SC‐ACOPF and standard ACOPF. (ii) Extend the GP‐based chance‐constraint approach to handle multiple simultaneous outages and more severe contingency scenarios, ensuring stability under $N$-$k$ disturbances. (iii) Leverage similar GP-based approaches to tackle network-wide small signal stability, dynamic stability of inverter-based resources, and other stability-affecting devices like FACTS and tap-changers. By pursuing these directions, we aim to further close the gap between rigorous stability assurance and computational efficiency, paving the way for secure, reliable operation of next‐generation power grids.

\begin{spacing}{1} 
\bibliography{references.bib}{}
\bibliographystyle{unsrt}
\end{spacing}

\end{document}